# On Tight Euclidean 6-Designs: An Experimental Result

**Djoko Suprijanto**

Combinatorial Mathematics Research Group
Faculty of Mathematics and Natural Sciences
Institut Teknologi Bandung 40132, INDONESIA.
Email: djoko@math.itb.ac.id

**Abstract.** A finite set $X \subseteq \mathbb{R}^n$ with a weight function $w: X \to \mathbb{R}_{>0}$ is called *Euclidean t-design* in $\mathbb{R}^n$ (supported by $p$ concentric spheres) if the following condition holds:

$$\sum_{i=1}^{p} \frac{w(X_i)}{|S_i|} \int_{S_i} f(x) d\sigma_i(x) = \sum_{x \in X} w(x) f(x)$$

for any polynomial $f(x) \in \text{Pol}(\mathbb{R}^n)$ of degree at most $t$. Here $S_i \subseteq \mathbb{R}^n$ is a sphere of radius $r_i \geq 0$, $X_i = X \cap S_i$, and $\sigma_i(x)$ is an $O(n)$-invariant measure on $S_i$ such that $|S_i| = r_i^{n-1} |S^{n-1}|$, with $|S_i|$ is the surface area of $S_i$ and $|S^{n-1}|$ is a surface area of the unit sphere in $\mathbb{R}^n$.

Recently, Bajnok (2006) [1] constructed tight Euclidean $t$-designs in the plane $(n=2)$ for arbitrary $t$ and $p$. In this paper we show that for case $t=6$ and $p=2$, tight Euclidean 6-designs constructed by Bajnok is the unique configuration in $\mathbb{R}^n$, for $2 \leq n \leq 8$.

## 1　　Introduction and Result

A combinatorial $t-(v,k,\lambda)$ design $X$ is one of important objects in combinatorics. It might be viewed (see [2]), in a sense, as an approximation of the discrete sphere $S_k$ of all *k*-subsets by the sub-collection $X$ of $S_k$, where

$$S_k := \{x \in \mathbb{R}^v : x_1^2 + x_2^2 + \cdots + x_v^2 = k, x_i \in \{0,1\}\}.$$

Later, Delsarte, Goethals and Seidel [3] introduced an analogue concept of designs for (continuous) sphere by defining what they called *spherical t-design*.

This new concept might be viewed as an approximation of the unit sphere $S^{n-1} \subseteq \mathbb{R}^n$ by the subset $X$ of $S^{n-1}$ with respect to integral of polynomial functions of degree at most $t$.

The concept of spherical $t$-design was generalized by Neumaier and Seidel ([4], see also Delsarte and Seidel [5]) by allowing weights and multiple spheres. In their papers, Neumaier and Seidel and also Delsarte and Seidel conjectured the non-existence of tight Euclidean $2e$-designs except the trivial ones.

**Conjecture 1.1** (Delsarte-Neumaier-Seidel). The only tight Euclidean $2e$-designs in $\mathbb{R}^n$, for $e \geq 2$, are regular simplices.

The first breakthrough on this area was performed by Bannai and Bannai [6]. Having slightly generalized the previous concept of Euclidean $t$-designs by dropping the condition of excluding **0** vector, they constructed a tight Euclidean 4-design in $\mathbb{R}^2$ supported by two concentric spheres as a counter-example for the conjecture. Moreover, they also completely classified tight Euclidean 4-designs with constant weight in $\mathbb{R}^n$, for $n \geq 2$, supported by two concentric spheres. Recently, in a joint work with Bannai and Bannai [7], the author introduced a new concept of strong non-rigidity for Euclidean $t$-designs. By using this new concept we also disproved Delsarte-Neumaier-Seidel's conjecture by showing the existence of infinitely many tight Euclidean designs having certain parameters.

Going back to Delsarte, Goethals and Seidel, regarding the spherical designs they showed that there is no tight spherical 6-designs in any Euclidean space except the one on $S^1$ (a unit sphere in a plane $\mathbb{R}^2$) ([3, Theorem 7.7]). Inspired by the situation in spherical designs, a natural question is what about tight Euclidean 6-designs? Are there any tight Euclidean 6-designs in $\mathbb{R}^n$? In answering this question, we divide them into two cases: the ones with constant weight and the others with non-constant weights. We observe that if the designs contain **0** vector, then by Lemma 2.15 given in the next section, $e = \frac{6}{2}$ should be even, which is impossible. Moreover, if $X$ is a tight Euclidean 6-design with constant weight, then Lemma 2.11(3), Remark 2.14, and Lemma 2.15 below imply $p = 2$ or 3. The purpose of this paper is to give a partial answer for the question by restricting our observation only to tight Euclidean 6-designs supported by two concentric spheres, sitting on the Euclidean spaces of small dimension $n$. Namely, we prove the following main theorem.

**Theorem 1.2** (Main theorem). *The only tight Euclidean 6-design in $\mathbb{R}^n$, for $2 \leq n \leq 8$, supported by two concentric spheres is the one in $\mathbb{R}^2$:*

$$X = \left\{ b_{kj} = \left( r_k \cos\left(\frac{2j+k}{5}\pi\right), r_k \sin\left(\frac{2j+k}{5}\pi\right) \right) : 1 \leq j \leq 5, 1 \leq k \leq 2 \right\}.$$

*The weight function of this design is $w(b_{kj}) = \frac{1}{r_k^5}$, for $k = 1, 2$.*

We remark that the design in the theorem above was constructed by Bajnok [1]. Hence the theorem says that in the Euclidean space of small dimension, Bajnok's construction of such designs is a unique configuration.

The paper is organized as follows. In section 2 we lay the groundwork for our result. We begin with some basic facts about association schemes. We also recall some facts about distance sets, spherical designs as well as Euclidean designs. We proved our main theorem in Section 3. Section 4 summarize the current status of classification of tight Euclidean designs. We end the paper by giving a conjecture on the (non-)existence of tight Euclidean 6-designs in $\mathbb{R}^n$, for $n \geq 2$, supported by two concentric spheres.

## 2     Preliminaries

This section contains some basic facts on association schemes, spherical designs, and Euclidean designs. We begin with association schemes.

### 2.1    Association Schemes

See [8] and [9] for undefined terms in association schemes.

Let $\mathfrak{X} = (X, \{R_i\}_{0 \leq i \leq d})$ be a symmetric association scheme. Let $P$ and $Q$ be the first and second eigenmatrix whose $(j,i)$-entry is $p_i(j)$ and $q_i(j)$, respectively. Denote $p_i(0) = k_i$ and $q_i(0) = m_i$ for $0 \leq i \leq d$. Then the following relation is well-known:

$$\frac{q_j(i)}{m_j} = \frac{p_i(j)}{k_j} \tag{2.1}$$

An association scheme $\mathfrak{X} = (X, \{R_i\}_{0 \leq i \leq d})$ is *imprimitive* if there exists a non-empty proper subset $\Lambda(\neq \{0\})$ of $\{0, 1, \ldots, d\}$ for which $\bigcup_{i \in \Lambda} R_i$ defines an equivalence relation on the set $X$. If $X$ is not imprimitive, the $\mathfrak{X}$ is called primitive.

**Lemma 2.1** (see, e.g., [10,11]). *Let* $\mathfrak{X} = \left( X, \{R_i\}_{0 \leq i \leq d} \right)$ *be a symmetric association scheme. For any i such that* $0 \leq i \leq d$, *the set* $X$ *can be embedded in the unit sphere* $S^{m_i - 1} \subseteq \mathbb{R}^{m_i}$, *where* $m_i = rank(E_i)$, *by*

$$\varphi : \begin{cases} X \to \mathbb{R}^{m_i} \\ x \mapsto \overline{x} := \sqrt{\dfrac{|X|}{m_i}} E_i e_x \end{cases}$$

*with* $e_x = (0, \ldots, 0, 1, 0, \ldots 0) \in \mathbb{R}^X$ *and 1 is in the x-th coordinate. If* $\mathfrak{X}$ *is primitive, then* $\varphi$ *is injective. Moreover,*

$$\langle \overline{x}, \overline{y} \rangle = \frac{q_i(j)}{m_i} = \frac{p_j(i)}{k_j}.$$

### 2.1.1 Krein Parameters

Since the Bose-Mesner algebra is also closed under the Hadamard product, then we may write

$$E_i \circ E_j = \frac{1}{|X|} \sum_{k=0}^{d} q_{ij}^k E_k$$

for some real numbers $q_{ij}^k$, called *Krein parameters*. These parameters are uniquely determined by the eigenmatrices $P$ and $Q$:

**Theorem 2.2** (Krein condition). *For all* $i, j, k \in \{0, 1, \ldots, d\}$ *we have,*

$$q_{ij}^k = \frac{1}{m_k |X|} \sum_{l=0}^{d} k_l q_i(l) q_j(l) q_k(l) \geq 0$$

where $k_l = p_l(0)$ and $m_k = q_k(0)$.

Let $\mathfrak{X} = \left( X, \{R_i\}_{0 \leq i \leq d} \right)$ be a symmetric association scheme. Let also $P$ and $Q$ be the eigenmatrices of the scheme. The association scheme $\mathfrak{X}$ is called *P*- polynomial [resp. $Q$-polynomial] scheme with respect to the ordering $R_0, R_1, \ldots R_d$ [resp. $E_0, E_1, \ldots E_d$] if there exist some polynomials $v_i(x)$ [resp. $v_i^*(x)$] of degree $i$ $(0 \leq i \leq d)$ such that $A_i = v_i(A_1)$ [resp. $E_i = v_i^*(E_1)$] under the Hadamard product]. Regarding the $Q$-polynomial scheme, Delsarte gives necessary and sufficient conditions for any symmetric association scheme to become a $Q$-polynomial scheme ([9, Theorem 5.16]), which together with the Krein condition can be restated as follows.

**Theorem 2.3** *A symmetric association scheme is $Q$-polynomial if and only if the Krein parameters $q_{ij}^k$ satisfy the following two conditions: $q_{1i}^{i+1} > 0$ and $q_{1i}^k = 0$ for $k > i+1$ $(0 \leq i \leq d-1)$.*

In particular, for the case $Q$-polynomial scheme of class 3, the following corollary is immediate.

**Corollary 2.4** *A symmetric association scheme is $Q$-polynomial of class 3 if and only if the Krein parameters $q_{ij}^k$ satisfy the conditions: $q_{11}^3 = 0$, $q_{11}^2 > 0$, and $q_{12}^3 > 0$.*

## 2.2   Distance Sets, Spherical and Euclidean Designs

### 2.2.1   Distance Sets

A finite subset $X \subseteq \mathbb{R}^n$ is called *antipodal* if $-x \in X$, for any $x \in X$. For a finite subset $X \subseteq \mathbb{R}^n$, we define

$$A(X) = \{\|x - y\| : x, y \in X, x \neq y\}.$$

$X$ is called an *s-distance* set if $|A(X)| = s$. For $\alpha \in A(X)$, we define

$$v_\alpha(x) = |\{y \in X : \|x - y\| = \alpha\}|.$$

Any subset $X$ is called *distance invariant* if $v_\alpha(x)$ does not depend on the choice of $x \in X$, and depend only on $\alpha$, for any fixed $\alpha \in A(X)$.

The upper bound for the cardinality of any $s$-distance set $X$ in $S^{n-1}$ is given by ([3, Theorem 4.8]):

$$|X| \leq \binom{n-1+s}{n-1} + \binom{n-2+s}{n-1}. \tag{2.2}$$

### 2.2.2  Spherical *t*-designs

Here is the exact definition of spherical t-designs as introduced by Delsarte, Goethals, and Seidel in 1977 [3].

**Definition 2.5**  Let *t* be a positive integer. A finite nonempty subset $X \subseteq S^{n-1}$ is called *spherical t-design* if the following condition holds:

$$\frac{1}{|S^{n-1}|} \int_{S^{n-1}} f(x) d\sigma(x) = \frac{1}{|X|} \sum_{x \in X} f(x) \tag{2.3}$$

for any polynomial $f(x) \in \mathbb{R}[x_1, x_2, \ldots, x_n]$ of degree at most t, where $\sigma(x)$ is the $O(n)$-invariant measure on $S^{n-1}$ and $|S^{n-1}|$ is the area of the sphere $S^{n-1}$. The maximum value of *t* for which $X$ is a spherical *t*-design is called the *strength* of $X$.

The lower bound for the cardinality of any spherical *t*-design $X$ in $S^{n-1}$ is given by ([3, Theorem 5.11 and Theorem 5.12]):

$$|X| \geq \binom{n-1+\left[\frac{t}{2}\right]}{n-1} + \binom{n-2+\left[\frac{t-1}{2}\right]}{n-1}. \tag{2.4}$$

An *s*-distance set [resp. spherical *t*-design] $X$ is called *tight* if the bound (2.2) [resp. (2.4)] is attained.

Again, the theorem below was also proved by Delsarte, Goethal, and Seidel [3].

**Theorem 2.6**  [3, Theorem 7.4] *Let $X$ be a spherical t-design as well as an s-distance set in $S^{n-1}$. If $t \geq 2s - 2$, then $X$ carries an s-class Q-polynomial scheme.*

**Remark 2.7**   In fact [3, Theorem 7.4], neither stated nor proved that the association scheme is *Q*-polynomial. The detail proof is given, e.g., in [12, Theorem 7.2.6], (c.f. [13, Theorem 9.6.4]) but Bannai and Bannai never claim the above theorem to be from them, instead they always refer the theorem to Delsarte, Goethal, and Seidel (see, e.g., [14, Theorem 3.5], [15]).

Now, let us turn to Euclidean designs.

### 2.2.3  Euclidean *t*-designs

Let $X$ be a finite set in $\mathbb{R}^n$, $n \geq 2$. Let $\{r_1, r_2, \ldots r_p\} = \{\|x\| : x \in X\}$, where $\|x\|$ is a norm of $x$ defined by standard inner product in $\mathbb{R}^n$ and $r_i$ is possibly 0. For each *i*, we define $S_i = \{x \in \mathbb{R}^n : \|x\| = r_i\}$, the sphere of radius $r_i$ centered at **0**. We say that $X$ is supported by the $p$ concentric spheres $S_1, S_2, \ldots, S_p$. If $r_i = 0$, then $S_i = \{\mathbf{0}\}$. Let $X_i = X \cap S_i$, for $1 \leq i \leq p$. Let $\sigma(x)$ be the $O(n)$-invariant measure on the unit sphere $S^{n-1} \subseteq \mathbb{R}^n$. We consider the measure $\sigma_i(x)$ on each $S_i$ so that $|S_i| = r_i^{n-1}|S^{n-1}|$, with $|S_i|$ is the surface area of $S_i$, namely $|S_i| = \dfrac{2 r_i^{n-1} \pi^{\frac{n}{2}}}{\Gamma\left(\frac{n}{2}\right)}$. We associate a positive real valued function $w$ on $X$, which is called *weight* of $X$. We define $w(X_i) = \sum_{x \in X} w(x)$. Here if $r_i = 0$, then we define $\frac{1}{|S_i|} \int_{S_i} f(x) d\sigma_i(x) = f(\mathbf{0})$, for any function $f(x)$ defined on $\mathbb{R}^n$. Let $S = \bigcup_{i=1}^p S_i$. Let $\varepsilon_S \in \{0,1\}$ be defined by

$$\varepsilon_S = \begin{cases} 1, & \mathbf{0} \in S \\ 0, & \mathbf{0} \notin S \end{cases}.$$

We give some more notation we use. Let $\mathrm{Pol}(\mathbb{R}^n) = \mathbb{R}[x_1, x_2, \ldots x_n]$ be the vector space of polynomials in *n* variables $x_1, x_2, \ldots x_n$. Let $\mathrm{Hom}_l(\mathbb{R}^n)$ be the subspace of $\mathrm{Pol}(\mathbb{R}^n)$ spanned by homogeneous polynomials of degree *l*. Let $\mathrm{Harm}_l(\mathbb{R}^n)$ be the subspace of $\mathrm{Pol}(\mathbb{R}^n)$ consisting of all harmonic polynomials. Let $\mathrm{Harm}_l(\mathbb{R}^n) = \mathrm{Harm}(\mathbb{R}^n) \cap \mathrm{Hom}_l(\mathbb{R}^n)$. Then we have $\mathrm{Pol}_l(\mathbb{R}^n) = \oplus_{i=0}^l \mathrm{Hom}_i(\mathbb{R}^n)$. Let $\mathrm{Pol}(S)$, $\mathrm{Pol}_l(S)$, $\mathrm{Hom}_l(S)$, $\mathrm{Harm}(S)$, and $\mathrm{Harm}_l(S)$ be the sets of corresponding polynomials restricted to the union $S$ of *p* concentric spheres. For example $\mathrm{Pol}(S) = \{f|_S : f \in \mathrm{Pol}(\mathbb{R}^n)\}$.

With the notation mentioned above, we define a Euclidean *t*-design as follows.

**Definition 2.8**  Let $X$ be a finite set with a weight function $w$ and let $t$ be a positive integer. Then $(X, w)$ is called *Euclidean t-design* in $\mathbb{R}^n$ if the following condition holds:

$$\sum_{i=1}^{p} \frac{w(X_i)}{|S_i|} \int_{S_i} f(x) d\sigma_i(x) = \sum_{x \in X} w(x) f(x),$$

for any polynomial $f(x) \in \text{Pol}(\mathbb{R}^n)$ of degree at most $t$. The maximum value of $t$ for which $X$ is a Euclidean $t$-design is called the *strength* of $X$.

The following theorem gives a condition which is equivalent to the definition of Euclidean $t$-designs.

**Theorem 2.9**  (Neumaier-Seidel). *Let $X$ be a finite nonempty subset in $\mathbb{R}^n$ with weight function w. Then the following (1) and (2) are equivalent:*

(1)  $X$ is a Euclidean t-design.
(2)  $\sum_{u \in X} w(u) \|u\|^{2j} \varphi(u) = 0$, *for any polynomial* $\varphi \in \text{Harm}_l(\mathbb{R}^n)$ *with* $1 \leq l \leq t$ *and* $0 \leq j \leq \left[\frac{t-l}{2}\right]$.

Let $X$ be a Euclidean $2e$-design in $\mathbb{R}^n$. Then it is known that ([4, Theorem 3.2], [5, Theorem 5.4]):

$$|X| \geq \dim(\text{Pol}_e(S)).$$

Following [6], we define the tightness for the Euclidean designs as given below.

**Definition 2.10**  Let $X$ be a Euclidean $2e$-design supported by $S$. If

$$|X| = \dim(\text{Pol}_e(S))$$

holds we call $X$ *tight Euclidean 2e-design on $S$*. Moreover, if

$$\dim(\text{Pol}_e(S)) = \dim(\text{Pol}_e(\mathbb{R}^n))$$

holds, then $X$ is called *tight Euclidean 2e-design*.

The following lemma is crucial in our study of Euclidean designs.

**Lemma 2.11**  [6, Lemma 1.10] *Let $X$ be a tight Euclidean 2e-design on $p$ concentric spheres. Then the following hold:*

(1) *The weight function $w$ is constant on each $X_i$, for $1 \leq i \leq p$.*

(2) *$X_i$ is at most an e-distance set $1 \leq i \leq p$.*

(3) *If the weight function $w$ is constant on $X \setminus \{0\}$, then $p - \varepsilon_S \leq e$.*

As an application of the above lemma, Et. Bannai [16] proved the following theorem. The theorem gives a certain connection between spherical designs and Euclidean designs.

**Theorem 2.12**  [16, Theorem 1.8] *Let $X$ be a tight Euclidean 2e-design on $p$ concentric spheres. If $p - \varepsilon_S \leq e$, then each $X_i$ is (similar to) a spherical $(2e - 2p + 2\varepsilon_S + 2)$-design. Moreover, if $p \leq \left[\frac{e + 2\varepsilon_S + 3}{2}\right]$, then each $X_i$ is a distance invariant set.*

Let $(X, w)$ be a finite weighted subset in $\mathbb{R}^n$. Let $S_1, S_2, \ldots S_p$ be the $p$ concentric spheres supporting $X$ and let $S = \bigcup_{i=1}^{p} S_i$.

For any $\varphi, \psi \in \text{Harm}(\mathbb{R}^n)$, we define the following inner-product

$$\langle \varphi, \psi \rangle = \frac{1}{|S^{n-1}|} \int_{S^{n-1}} \varphi(x) \psi(x) d\sigma(x).$$

Then we have the following (see [17], [3],[5],[12],[16]).

**Lemma 2.13**  The following three statements hold:

(1) $\text{Harm}(\mathbb{R}^n)$ is a positive definite inner-product space under $\langle -, - \rangle$ and has the orthogonal decomposition $\text{Harm}(\mathbb{R}^n) = \perp_{i=0}^{\infty} \text{Harm}_i(\mathbb{R}^n)$.

(2) $\text{Pol}_e(\mathbb{R}^n) = \bigoplus_{0 \leq i + 2j \leq e} \|x\|^{2j} \text{Harm}_i(\mathbb{R}^n)$ with $\dim(\text{Pol}_e(\mathbb{R}^n)) = \binom{n+e}{e}$.

(3)

$$\mathrm{Pol}_e(S) = \left\langle \|x\|^{2j} : 0 \le j \le \min\left\{p-1, \left[\tfrac{e}{2}\right]\right\}\right\rangle \oplus \left\{ \bigoplus_{\substack{1 \le i \le e \\ 0 \le j \le \min\left\{p-\varepsilon_S-1, \left[\tfrac{e-1}{2}\right]\right\}}} \|x\|^{2j} \mathrm{Harm}_i(S) \right\}$$

.

(a) If $p \le \left[\tfrac{e+\varepsilon_S}{2}\right]$ then

$$\dim(\mathrm{Pol}_e(S)) = \varepsilon_S + \sum_{i=0}^{2(p-\varepsilon_S)-1} \binom{n+e-i-1}{n-1} < \sum_{i=0}^{e} \binom{n+e-i-1}{n-1} = \binom{n+e}{e},$$

where $e$ is a non-negative integer.

(b) If $p \ge \left[\tfrac{e+\varepsilon_S}{2}\right]+1$, then

$$\dim(\mathrm{Pol}_e(S)) = \binom{n+e}{e},$$

where $e$ is a non-negative integer.

**Remark 2.14** Definition 2.10 and Lemma 2.13 show that a tight Euclidean $2e$-design is the same as a tight Euclidean $2e$-design on $p$ concentric spheres with $p \ge \left[\tfrac{e+\varepsilon_S}{2}\right]+1$.

The next lemma is stated in Bannai and Bannai [6].

**Lemma 2.15** [6, Proposition 1.7] *Let $X \in \mathbb{R}^n$ be a tight Euclidean $2e$-design. If $0 \in X$, then $e$ is even, $p = \tfrac{e}{2}+1$, and $X \setminus \{0\}$ is a tight Euclidean $2e$-design on $\tfrac{e}{2}$ concentric spheres.*

Let $h_l = \dim(\mathrm{Harm}_l(\mathbb{R}^n))$ and $\varphi_{l,1}, \ldots, \varphi_{l,h_l}$ be an orthonormal basis of $\mathrm{Harm}_l(\mathbb{R}^n)$ with respect to the inner-product defined above. Then, by Lemma 2.13,

$$\left\{\|x\|^{2j} : 0 \le j \le \min\left\{p-1, \left[\tfrac{e}{2}\right]\right\}\right\} \cup$$

$$\left\{\|x\|^{2j} \varphi_{l,i}(x) : 1 \le l \le e, 1 \le i \le h_l, 0 \le j \le \min\left\{p-\varepsilon_S-1, \left[\tfrac{e-l}{2}\right]\right\}\right\}$$

gives a basis of $\text{Pol}_e(S)$.

Now, we are going to construct a more convenient basis of $\text{Pol}_e(S)$ for our purpose. Let $\mathcal{G}(\mathbb{R}^n)$ be the subspace of $\text{Pol}_e(S)$ spanned by $\{\|x\|^{2j}:0\leq j\leq p-1\}$. Let $\mathcal{G}(X)=\{g\,|_X: g\in\mathcal{G}(\mathbb{R}^n)\}$. Then $\{\|x\|^{2j}:0\leq j\leq p-1\}$ is a basis of $\mathcal{G}(X)$. We define an inner-product $\langle -,-\rangle_l$ on $\mathcal{G}(X)$ by

$$\langle f,g\rangle_l = \sum_{x\in X} w(x)\|x\|^{2l} f(x)g(x), \text{ for } 1\leq l\leq e. \tag{2.5}$$

We apply the Gram-Schmidt method to the basis $\{\|x\|^{2j}:0\leq j\leq p-1\}$ to construct an orthonormal basis

$$\{g_{l,0}(x), g_{l,1}(x),\ldots, g_{l,p-1}(x)\}$$

of $\mathcal{G}(X)$ with respect to the inner-product $\langle -,-\rangle_l$. We can construct them so that for any $l$ the following holds:

$g_{l,j}(x)$ is a linear combination of $1, \|x\|^2,\ldots,\|x\|^{2j}$, with $\deg(g_{l,j})=2j$ for $0\leq j\leq p-1$.

As an example, for $p=2$, we can express $g_{l,j}$ in the following way:

$$g_{l,0}(x) \equiv \frac{1}{\sqrt{a_l}},\ g_{l,1} = \sqrt{\frac{a_l}{a_l a_{l+2}-a_{l+1}^2}}\left(\|x\|^2 - \frac{a_{l+1}}{a_l}\right) \tag{2.6}$$

with $a_l = \sum_{x\in X} w(x)\|x\|^{2l}$.

Now we are ready to give a new basis for $\text{Pol}_e(S)$. Let us consider the following sets:

$$\mathcal{H}_0 = \left\{ g_{0,j} : 0 \leq j \leq \min\left\{ p-1, \left[\frac{e}{2}\right] \right\} \right\},$$

$$\mathcal{H}_l = \left\{ g_{l,j} \varphi_{l,i} : 0 \leq j \leq \min\left\{ p - \varepsilon_S - 1, \left[\frac{e-l}{2}\right] \right\}, 1 \leq i \leq h_l \right\}, \text{ for } 1 \leq l \leq e.$$

Then $\mathcal{H} = \bigcup_{l=0}^{e} \mathcal{H}_l$ is a basis of $\text{Pol}_e(S)$.

We close this section by the following lemma (see [6] or [7] for the proof.)

**Lemma 2.16.** *If* $(X, w)$ *is a tight Euclidean 2e-design on* $S$, *then the following* (1) *and* (2) *hold:*

(1) *The weight function of* $X$ *satisfies*

$$\sum_{\substack{1 \leq l \leq e, \\ 0 \leq j \leq \min\{p-\varepsilon_S-1, [\frac{e-l}{2}]\}}} \|u\|^{2l} g_{l,j}^2(u) Q_l(1) + \sum_{j=0}^{\min\{p-1, [\frac{e}{2}]\}} g_{0,j}^2(u) = \frac{1}{w(u)}, \text{ for all } u \in X. \quad (2.7)$$

(2) *For any distinct points* $u, v \in X$, *we have*

$$\sum_{\substack{1 \leq l \leq e, \\ 0 \leq j \leq \min\{p-\varepsilon_S-1, [\frac{e-l}{2}]\}}} \|u\|^l \|v\|^l g_{l,j}(u) g_{l,j}(v) Q_l\left(\frac{\langle u,v \rangle}{\|u\|\|v\|}\right) + \sum_{j=0}^{\min\{p-1, [\frac{e}{2}]\}} g_{0,j}(u) g_{0,j}(v) = 0.$$

*(2.8)*

*Here* $\langle u, v \rangle$ *is the standard inner-product in Euclidean space* $\mathbb{R}^n$ *and* $Q_l(\alpha)$ *is the Gegenbauer polynomial of degree l.*

## 3 Proof of Main Theorem

We prove Theorem 1.2 by contradiction. The general idea is to show that the assumption of the existence of tight Euclidean 6-design of certain given parameters does not carry a *Q*-polynomial scheme of class 3. Hence we get a contradiction. The detail follows.

Let $X = X_1 \cup X_2$ be a tight Euclidean 6-design in $\mathbb{R}^n$, for $2 \leq n \leq 8$. By Theorem 2.12, we know that $X_i$ is (similar to) a spherical 4-design. We also know, by Lemma 2.11, that a tight spherical 4-design $X_i$ is also a 2-distance set, while the non-tight one is also a 3-distance set. Therefore, Theorem 2.6 guarantees that the non-tight spherical 4-design $X_i$ should carry a 3-class $Q$-polynomial scheme.

On the other hand, van Dam [18] gives all feasible character tables of the 3-class symmetric association schemes on points up to 100. By the help of Lemma 2.1, we know that the symmetric association scheme on $X_i$ can be embedded into a unit sphere, which also give us the feasible 3-inner product set. Hence, keeping in mind that any distance set in the unit sphere has one-to-one correspondence with an inner product set, we can investigate whether the finite set $X_i$ carries a 3-class $Q$-polynomial scheme, by comparing the numerical 3-inner product sets (obtained from Lemma 2.16) with the feasible ones (given by van Dam's character tables).

Let us consider first some special cases.

## 3.1 Some Special Cases

We begin with some elementary facts. Let $N, N_i$ denotes the cardinality of $X, X_i$, for $i = 1, 2$, respectively. Suppose $N_1 \leq N_2$. Then the lower bound $L_b$ and upper bound $U_b$ of $N_1$ is given by

| $n$ | 2 | 3 | 4 | 5 | 6 | 7 | 8 |
|---|---|---|---|---|---|---|---|
| $N$ | 10 | 20 | 35 | 56 | 84 | 120 | 165 |
| $[L_b, U_b]$ | 5 | [9,10] | [14,17] | [20,28] | [27,42] | [35,60] | [44,82] |

We notice that there are three kind 3-class symmetric association schemes of "degenerate" case (to follow van Dam [vDam-99]). They are: (1) the schemes generated by $n$ disjoint union of strongly regular graphs $SRG(v, k, \lambda, \mu)$ (2) the schemes generated by $SRG(v, k, \lambda, \mu) \otimes J_n$, and (3) the rectangular scheme $R(m, n)$. The character tables of each case are (see [18, p. 88]):

$$(1) \begin{pmatrix} 1 & k & v-1-k & (n-1)v \\ 1 & k & v-1-k & -v \\ 1 & r & -1-r & 0 \\ 1 & s & -1-s & 0 \end{pmatrix}, (2) \begin{pmatrix} 1 & nk & n-1 & n(v-1-k) \\ 1 & nr & n-1 & n(-1-r) \\ 1 & 0 & -1 & 0 \\ 1 & ns & n-1 & n(-1-s) \end{pmatrix}$$

and

$$(3) \begin{pmatrix} 1 & (m-1)(n-1) & n-1 & m-1 \\ 1 & 1 & -1 & -1 \\ 1 & 1-m & -1 & m-1 \\ 1 & 1-n & n-1 & -1 \end{pmatrix}.$$

By direct calculation it is easy to see that some values of Krein parameters are

$$(1)\ q_{11}^3 = 0, q_{11}^2 = 0,\ (2)\ q_{11}^2 = 0, q_{12}^3 = 0$$

$$(3)\ q_{11}^3 = (n-1)(m-2), q_{11}^2 = (n-2)(m-1), q_{12}^3 = n-1.$$

Hence for the first two cases, the schemes are not $Q$-polynomial, while for the last case, the scheme is $Q$-polynomial if and only if $n = 2$ or $m = 2$, that is if cardinality of the finite set carrying the scheme is even. Furthermore, the only feasible 3-inner product set given by this scheme is $\{\pm\frac{1}{m-1}, -1\}$, for $m \neq 2$. We will include this feasible set in our observation below.

**Remark 3.1** We checked the above degenerate cases for all possible ordering of the primitive idempotent basis $E_0$, $E_1$, $E_2$, $E_3$. Here we consider the character table $P'$ obtained from $P$ by applying a permutation to the set of its rows but the first (there are six possibilities) and we have: two of them give $q_{11}^3 = 0$, $q_{11}^2 = 0$ and the others $q_{11}^3 \neq 0$. (for type (1) and (2)); and two of them give $q_{11}^3 = (n-1)(m-2)$, $q_{11}^2 = (n-2)(m-1)$, $q_{12}^3 = n-1$ and the others $q_{11}^3 = 0$, $q_{11}^2 = 0$ (for type (3)).

Next, let us consider the following cases.

### 3.1.1  Case $(n, N_1) = (3, 9), (4, 14), (5, 20), (7, 35), (8, 44)$.

These are the cases where $X_1$ is a tight spherical 4-design. By Bannai-Damarell's criteria (see [19], c.f. [6, Remark 4]), it is known that if a tight spherical 4-design $X_1 \subseteq S^{n-1}$ exists, then $n = (2m+1)^2 - 3$ holds, for some integer $m$. Since there is no $m$ satisfying $(2m+1)^2 = 6$, [resp. 7, 8, 10, and 11], then there does not exist tight spherical 4-design on $S^2$, [resp. $S^3$, $S^4$, $S^6$, and $S^7$].

Hence there is no tight Euclidean 6-design in $\mathbb{R}^3$ [resp. $\mathbb{R}^4$, $\mathbb{R}^5$, $\mathbb{R}^7$, and $\mathbb{R}^8$] supported by two concentric spheres with such parameters.

**Remark 3.2** Results for the cases $(n, N_1) = (3, 9)$, $(4, 14)$, and $(5, 20)$ are also a direct consequence of a work of Boyvalenkov and Nikova [20], where they improved the lower bound of tight spherical 4-designs on $S^2$, $S^3$, and $S^4$, from 9, 14, and 20 to 10, 15, and 21, respectively.

Now, let us turn to the general treatment. We begin with the constant weight case.

### 3.2    Case 1: Constant Weight

For $n = 2$, then $N_i = 5$, i.e., $X_i$ are regular pentagon (5-gon). We may assume that $\|x\| = 1$, for $x \in X_1$ and $w(x) = 1$, for $x \in X$ (see [6, Proposition 2.4], c.f. Theorem 2.9). From equation (2.7) we have that the weight function for any point sitting on the second sphere is $w(x) = \frac{1}{\sqrt{\|x\|^5}}$. Our assumption implies $\|x\| = 1$, for any $x \in X_2$, which is impossible. Hence tight Euclidean 6-design in $\mathbb{R}^2$ supported by two concentric spheres does not exist.

Next, let us consider the other cases. We mention first the procedure we use in our observation. Let $X = X_1 \cup X_2$ be a tight Euclidean 6-design with constant weight in $\mathbb{R}^n$ $(2 \leq n \leq 8)$, let $|X_1| = N_1$, and $|X_2| = N_2$. We may assume that $\|x\| = 1$, for $x \in X_1$, and $w(x) = 1$, for $x \in X$.

**Step 1:** Given $n$, $N_1$, and $N_2$.
**Step 2:** If there exist 3-class symmetric association schemes on $N_1$ and $N_2$ points simultaneously, then further check if $X_i$ carries a $Q$-polynomial scheme (by Corollary 2.4). If such a polynomial scheme exists then
– Calculate the radius of the second sphere $\|x\| = \sqrt{R}$, for $x \in X_2$ by equation (2.7) in Lemma 2.16.
– If $R \neq 1$, then substitute $R$ to equation (2.8) in Lemma 2.16 to get the 3-inner product set $\{\alpha_1, \alpha_2, \alpha_3\}$.
**Step 3:** Compare $\{\alpha_1, \alpha_2, \alpha_3\}$ with the entries of corresponding character table $\left\{\frac{p_1(i)}{k_1}, \frac{p_2(i)}{k_2}, \frac{p_3(i)}{k_3}\right\}$, for $2 \leq i \leq 4$.

(See the Appendix for example of calculation results.)

### 3.3    Case 2: Non-constant Weight

For *n*=2, then $N_i = 5$, i.e. $X_i$ are regular 5-gons. Bajnok [1] has constructed an example of such a design:

$$X = \left\{ b_{kj} = \left( r_k \cos\left( \frac{2j+k}{5}\pi \right), r_k\left( \frac{2j+k}{5}\pi \right) \right) : 1 \leq j \leq 5, 1 \leq k \leq 2 \right\}.$$

The weight function is given by $w(b_{kj}) = \frac{1}{r_k^5}$, for $k = 1,2$. Hence there exists a tight Euclidean 6-design in $\mathbb{R}^2$ supported by two concentric spheres. Moreover, it is easy to show that in fact it is the only tight Euclidean 6-design in $\mathbb{R}^2$. The argument is as follows. Let $X' = X'_1 \cup X'_2$ be a tight Euclidean 6-design in $\mathbb{R}^2$. Then $X'_1$ and $X'_2$ should be tight spherical 4-designs, namely regular pentagons. Hence, up to the action of orthogonal group $O(2)$, we can write $X'_1 = X_1$ and $X'_2 = \varphi(X_2)$, for some rotation $\varphi \in O(2)$. Using the Neumaier-Seidel's Theorem above it can be shown that $\varphi = I$, the identity, namely $X'_2 = X_2$.

Next, let us consider the other cases. We mention first the procedure we use in our observation. Let $X = X_1 \cup X_2$ be a tight Euclidean 6-design with non-constant weight in $\mathbb{R}^n (3 \leq n \leq 8)$, let $|X_1| = N_1$ and $|X_2| = N_2$. Again, we may assume that $\|x\| = 1$, for $x \in X_1$ and $w(x) = 1$, for $x \in X_1$.

**Step 1:** Given *n*, $N_1$, and $N_2$.
**Step 2:** If there exist 3-class symmetric association schemes on $N_1$ and $N_2$ points simultaneously, then further check if $X_i$ carries a *Q*-polynomial scheme (by Corollary 2.4). If such a polynomial scheme exists, then
   – Calculate the weight function *w* by equation (7) in Lemma 2.16.
   – Substitute the weight function *w* to equation (8) in Lemma 2.16 to get the 3-inner product set $\{\alpha_1, \alpha_2, \alpha_3\}$. (Here $\alpha_i (1 \leq i \leq 3)$ are functions of the second radius *R*).
**Step 3:** Compare $\{\alpha_1, \alpha_2, \alpha_3\}$ with the entries of corresponding character table $\left\{ \frac{p_1(i)}{k_1}, \frac{p_2(i)}{k_2}, \frac{p_3(i)}{k_3} \right\}$, for $2 \leq i \leq 4$. (Here we consider at most $3 \times 9 = 27$ equations for one corresponding character table).

- If, say, $\alpha_1 = \frac{p_1(i)}{k_1}$ implies the positive real value of radius $R_0 (\neq 1)$, then substitute $R_0$ to $\alpha_2$ and $\alpha_3$.
- Check whether $\{\alpha_1(R_0), \alpha_2(R_0), \alpha_3(R_0)\} = \left\{\frac{p_1(i)}{k_1}, \frac{p_2(i)}{k_2}, \frac{p_3(i)}{k_3}\right\}$, for some $i \in \{2, 3, 4\}$.

(See the Appendix for example of calculation results.)

In summary, our assumption of the existence of tight Euclidean 6-designs with certain given parameters implies:

1. the non-existence of tight spherical 4-design in a Euclidean space of given dimension (by Bannai-Damarell's criteria), or
2. the non-existence of 3-class symmetric association scheme on $N_1$ or $N_2$ points (by checking on van Dam's table), or
3. the non-existence of 3-class $Q$-polynomial scheme on $N_1$ or $N_2$ points (by Krein condition), or
4. the 3-class symmetric association scheme on $N_1$ points does not provide the 3-inner product set (by looking at van Dam's table), or
5. the numerical 3-inner product set does not appear, or
6. the numerical 3-inner product set does not coincide with the 3-inner product set provided by the character table of $Q$-polynomial scheme on $N_1$ points.

All of these lead to a contradiction. Hence, we have proved the main theorem.

## 4 Concluding Remarks

As we have seen, there is no tight Euclidean 6-design supported by two concentric spheres in Euclidean spaces of small dimensions, namely in $\mathbb{R}^n$ $(2 \leq n \leq 8)$, for almost all feasible parameters. The only exception is tight Euclidean 6-designs of non-constant weight in $\mathbb{R}^2$. The designs was constructed by Bajnok ([1, Theorem 9]). Our effort here might be regarded as a continuation of the work on giving classifications of tight Euclidean designs in $\mathbb{R}^n$. Hence, the current status of this work is as follows:

1. tight Euclidean 2-designs in $\mathbb{R}^n$, for $n \geq 2$, supported by all feasible concentric spheres [7].
2. tight Euclidean 3-designs in $\mathbb{R}^n$, for $n \geq 2$, supported by all feasible concentric spheres [16].

3. tight Euclidean 4-designs in $\mathbb{R}^n$, for $n \geq 2$, with constant weight, supported by two concentric spheres [6].
4. tight Gaussian 4-designs, i.e. a special kind of tight Euclidean designs, in $\mathbb{R}^n$, for $n \geq 2$ supported by two concentric spheres [21].
5. tight Euclidean 5-designs in $\mathbb{R}^n$, for $n \geq 2$ supported by two concentric spheres [16].
6. tight Euclidean 6-designs in $\mathbb{R}^n$, for $2 \leq n \leq 8$, supported by two concentric spheres.

Recently Et. Bannai [14] has constructed some examples of (non-constant weight) tight Euclidean 4-designs supported by two concentric spheres, but the problem of classification of such designs still far of being solved. Besides, Bajnok [1, 22] has also constructed some other sporadic examples of Euclidean designs. Very recently, Bannai, Bannai, Hirao, and Sawa [23] announced that they have classified tight Euclidean 7-designs supported by two concentric spheres completely, but the complete proof is now in preparation by Bannai and Bannai.

Regarding the tight Euclidean 6-designs on two concentric spheres, all existing phenomena mentioned above lead us to a conjecture on the complete classification of such designs. We end this paper by the conjecture:

**Conjecture 4.1** The only tight Euclidean 6-designs in $\mathbb{R}^n$, for $n \geq 2$ supported by two concentric spheres are:

$$X = \left\{ b_{kj} = \left( r_k \cos\left( \frac{2j+k}{5}\pi \right), r_k \sin\left( \frac{2j+k}{5}\pi \right) \right) : 1 \leq j \leq 5, 1 \leq k \leq 2 \right\}.$$

The weight function of these designs are $w(b_{kj}) = \frac{1}{r_k^5}$ for $k = 1, 2$.

**Remark 4.2** One of the referee, Eiichi Bannai, sent me the paper [11], where they give a new example of tight Euclidean 6-design in $\mathbb{R}^{22}$. It becomes the first example of tight Euclidean 6-design in $n$-dimensional Euclidean space, for $3 \leq n \leq 438$, supported by two concentric spheres. It also gives a counter-example of the above Conjecture. However, Bannai et.al. seem to believe that it is the only existing tight Euclidean 6-design in $\mathbb{R}^n$, for $n \geq 3$. Hence, we modify the conjecture to the following question.

**Question 4.2.** Find another example of tight Euclidean 6-design in $\mathbb{R}^n$, for $n \geq 3$, on two concentric spheres, or prove that there exist no such tight Euclidean 6-design except the ones given in the above conjecture and in [11].

## Acknowledgement

The author would like to thank Prof. Eiichi Bannai and Prof. Etsuko Bannai for introducing the problem and many valuable suggestions. The author would also like to thank the referees for careful reading of the manuscript and many helpful comments.

**Appendix: Example of Calculation Results**

The table below contains our calculation results. It shows that in almost all cases, our assumption on the existence of tight Euclidean 6-designs having given parameters do not imply the structure of $Q$-polynomial schemes for $X_1$ or $X_2$. In some cases, either the tight Euclidean 6-designs in certain given parameters do not give 3-inner product set, or provide numerical 3-inner product set which does not coincide with the feasible 3-inner product set. All of these results lead to a contradiction.

In the following tables, "feasible" 3-inner product sets refer to the sets which are coming from the first eigenmatrix $P$ of related 3-class symmetric association schemes, as given by van Dam [18], while "numerical" 3-inner product sets come from the calculation results. The number $[\alpha]$ in the "feasible" column refers to the page where the sets are mentioned in [18]. $\times_\alpha$ in the "Krein parameters" column means the 3-class symmetric association scheme on $\alpha$ points does not exist. We give here only the case $n = 7$, for only one ordering of primitive idempotent basis $E_0, E_1, E_2, E_3$. The treatment for another cases is exactly the same.

**Case $n = 7$**

Here we have $N := 120$ and $35 \leq N_1 \leq 60$. We have shown that there exist no tight Euclidean 6-design in $\mathbb{R}^7$ supported by two concentric spheres, consisting of $N_1 = 35$ (see Case 3.1.1 for the detail).

| $N_1$ | Krein parameters | feasible 3-i.p.set |
|---|---|---|
| 43 | $\times_{77}$ | |
| 44 | $q_{11}^3 = 0.000, q_{11}^2 = 0.000, q_{12}^3 = 1.000$ | $\left\{\frac{2.236}{7}, -\frac{1}{21}, -\frac{2.236}{15}\right\}$ [100], $\left\{-\frac{2.236}{7}, -\frac{1}{21}, \frac{2.236}{15}\right\}$ [100] |
| | $q_{11}^3 = 0, q_{11}^2 = 0, q_{12}^3 = 11.00$ | $\left\{\frac{4.583}{21}, -\frac{4.583}{21}, -1\right\}$ [100] |
| 45 | $q_{11}^3 = \frac{325}{32}, q_{11}^2 = \frac{425}{32}, q_{12}^3 = \frac{225}{32}$ | $\left\{\frac{1}{4}, -\frac{1}{20}, -\frac{1}{8}\right\}$ [93], $\left\{-\frac{1}{2}, \frac{1}{4}, -\frac{1}{8}\right\}$ [93] |
| | $q_{11}^3 = \frac{3}{2}, q_{11}^2 = 3, q_{12}^3 = 3$ | $\left\{\frac{1}{2}, -\frac{1}{8}, -\frac{1}{4}\right\}$ [93] |
| | $q_{11}^3 = \frac{105}{32}, q_{11}^2 = \frac{165}{32}, q_{12}^3 = \frac{75}{8}$ | $\left\{\frac{1}{4}, -\frac{1}{8}, -\frac{1}{4}\right\}$ [93] |
| | $q_{11}^3 = \frac{15}{2}, q_{11}^2 = 10, q_{12}^3 = \frac{25}{2}$ | $\left\{\frac{1}{8}, -\frac{1}{8}, -\frac{1}{4}\right\}$ [93] |
| | $q_{11}^3 = \frac{513}{32}, q_{11}^2 = \frac{621}{32}, q_{12}^3 = \frac{27}{4}$ | $\left\{\frac{1}{12}, -\frac{1}{8}, -\frac{1}{4}\right\}$ [93] |
| | $q_{11}^3 = 7.500, q_{11}^2 = 7.500, q_{12}^3 = 7.500$ | $\left\{\frac{4.873}{16}, -\frac{2.873}{16}, -\frac{1}{4}\right\}$ [100] |
| 46 | $q_{11}^3 = 0.000, q_{11}^2 = 0.000, q_{12}^3 = 1.000$ | $\left\{\frac{2.449}{11}, -\frac{1}{22}, -\frac{2.449}{12}\right\}$ [100], $\left\{-\frac{2.449}{11}, -\frac{1}{22}, \frac{2.449}{12}\right\}$ [100] |
| 47 | $\times_{47}$ | |
| 48 | $q_{11}^3 = 15, q_{11}^2 = \frac{52}{3}, q_{12}^3 = 15$ | $\left\{\frac{1}{6}, -\frac{1}{15}, -\frac{1}{10}\right\}$ [93], $\left\{-\frac{1}{3}, -\frac{1}{15}, \frac{1}{5}\right\}$ [93] |
| | $q_{11}^3 = \frac{3}{2}, q_{11}^2 = 4, q_{12}^3 = 6$ | $\left\{\frac{1}{3}, -\frac{1}{6}, -\frac{1}{2}\right\}$ [93], $\left\{-\frac{1}{5}, \frac{1}{10}, -\frac{1}{2}\right\}$ [93] |
| 49 | $\times_{71}$ | |

| $N_1$ | Krein parameters | feasible 3-i.p.set |
|---|---|---|
| 50 | $q_{11}^3 = 0.000, q_{11}^2 = 0.000, q_{12}^3 = 1.000$ | $\left\{\frac{2.449}{9}, -\frac{1}{24}, -\frac{2.449}{16}\right\}$ [100], $\left\{-\frac{2.449}{9}, -\frac{1}{24}, \frac{2.449}{16}\right\}$ [100] |
| 51 | $q_{11}^3 = \frac{51}{16}, q_{11}^2 = \frac{255}{32}, q_{12}^3 = \frac{17}{2}$ | $\left\{\frac{1}{4}, -\frac{1}{8}, -\frac{1}{2}\right\}$ [93], $\left\{-\frac{1}{4}, \frac{1}{8}, -\frac{1}{2}\right\}$ [93] |
| 52 | $q_{11}^3 = 0, q_{11}^2 = \frac{156}{25}, q_{12}^3 = 13 (*)$ | $\left\{\frac{1}{5}, -\frac{1}{5}, -1\right\}$ [93] |
|  | $q_{11}^3 = 13.42, q_{11}^2 = 19.81, q_{12}^3 = 6.000$ | $\left\{-\frac{1}{3}, \frac{1}{9}, -\frac{1}{27}\right\}$ [100], $\left\{\frac{3.732}{6}, \frac{0.464}{18}, -\frac{5.196}{27}\right\}$ [100], $\left\{\frac{0.268}{6}, -\frac{6.464}{18}, \frac{5.196}{27}\right\}$ [100] |
|  | $q_{11}^3 = 4.000, q_{11}^2 = 6.000, q_{12}^3 = 7.000$ | $\left\{\frac{4.302}{17}, -\frac{1.548}{17}, -\frac{3.754}{17}\right\}$ [104] |
| 53 | $\times_{53}$ | |
| 54 | $q_{11}^3 = 0.000, q_{11}^2 = 0.000, q_{12}^3 = 10.00$ | $\left\{\frac{2.646}{13}, -\frac{1}{26}, -\frac{2.646}{14}\right\}$ [100], $\left\{-\frac{2.646}{13}, -\frac{1}{26}, \frac{2.646}{14}\right\}$ [100] |
| 55 | $q_{11}^3 = 1.358, q_{11}^2 = 1.358, q_{12}^3 = 5.432$ | $\left\{\frac{3.854}{18}, -\frac{1}{9}, -\frac{2.854}{18}\right\}$ [100] |
| 56 | $q_{11}^3 = 0, q_{11}^2 = \frac{448}{225}, q_{12}^3 = \frac{14}{5} (*)$ | $\left\{\frac{7}{15}, -\frac{1}{15}, -\frac{3}{5}\right\}$ [93], $\left\{\frac{1}{15}, -\frac{1}{6}, \frac{3}{10}\right\}$ [93], $\left\{-\frac{1}{5}, \frac{1}{10}, -\frac{1}{10}\right\}$ [93] |
|  | $q_{11}^3 = 0, q_{11}^2 = \frac{140}{9}, q_{12}^3 = 21(*)$ | $\left\{\frac{1}{9}, -\frac{1}{9}, -1\right\}$ [94], $\left\{-\frac{1}{3}, \frac{1}{3}, -1\right\}$ [94] |
|  | $q_{11}^3 = 2.718, q_{11}^2 = 7.568, q_{12}^3 = 6.000$ | $\left\{-\frac{3}{5}, \frac{1}{5}, -\frac{1}{15}\right\}$ [100], $\left\{\frac{2.414}{5}, \frac{0.828}{20}, -\frac{4.243}{30}\right\}$ [100], $\left\{-\frac{0.414}{5}, -\frac{4.828}{20}, \frac{4.243}{30}\right\}$ [100] |
| 57 | $q_{11}^3 = 7.355, q_{11}^2 = 11.59, q_{12}^3 = 6.000$ | $\left\{-\frac{1}{2}, \frac{1}{10}, -\frac{1}{20}\right\}$ [100], $\left\{\frac{2.618}{6}, \frac{0.854}{30}, -\frac{4.472}{20}\right\}$ [100], $\left\{\frac{0.382}{6}, -\frac{5.854}{30}, \frac{4.472}{20}\right\}$ [100] |
| 58 | $q_{11}^3 = 0.000, q_{11}^2 = 0.000, q_{12}^3 = 1.000$ | $\left\{\frac{2.449}{8}, -\frac{1}{28}, -\frac{2.449}{21}\right\}$ [100], $\left\{-\frac{2.449}{8}, -\frac{1}{28}, \frac{2.449}{21}\right\}$ [100] |

| $N_1$ | Krein parameters | feasible 3-i.p.set |
|---|---|---|
| 59 | $\times_{59}$ | |
| 60 | $q_{11}^3 = \frac{25}{3}, q_{11}^2 = \frac{25}{2}, q_{12}^3 = \frac{20}{3}$ | $\left\{\frac{1}{5}, 0, -\frac{1}{5}\right\}$ [94], $\left\{0, -\frac{1}{4}, \frac{1}{4}\right\}$ [94], $\left\{-\frac{1}{3}, \frac{1}{6}, 0\right\}$ [94] |
| | $q_{11}^3 = \frac{96}{7}, q_{11}^2 = \frac{816}{49}, q_{12}^3 = \frac{128}{7}$ | $\left\{\frac{1}{7}, -\frac{1}{14}, -\frac{1}{8}\right\}$ [94], $\left\{-\frac{4}{21}, -\frac{1}{14}, \frac{1}{6}\right\}$ [94] |
| | $q_{11}^3 = 0, q_{11}^2 = 0, q_{12}^3 = 6.000$ | $\left\{\frac{3.317}{11}, -\frac{3.317}{44}, -\frac{1}{4}\right\}$ [100], $\left\{-\frac{3.317}{11}, \frac{3.317}{44}, -\frac{1}{4}\right\}$ [100] |
| | $q_{11}^3 = 0, q_{11}^2 = 0, q_{12}^3 = 10.00$ | $\left\{\frac{4.359}{19}, -\frac{4.359}{38}, -\frac{1}{2}\right\}$ [100], $\left\{-\frac{4.359}{19}, \frac{4.359}{38}, -\frac{1}{2}\right\}$ [100] |
| | $q_{11}^3 = 0, q_{11}^2 = 0, q_{12}^3 = 15$ | $\left\{\frac{5.385}{29}, -\frac{5.385}{29}, -1\right\}$ [100] |

We continue to step 3, for some tight Euclidean 6-designs carrying $Q$-polynomial schemes (marked by $(*)$), by comparing feasible and numerical 3-inner product set. Consider the table below.

| [$N, N_1$] | Feasible | 3-inner product set | |
|---|---|---|---|
| | | Numerical | |
| | | constant weight | non-constant weight |
| [120, 36] | $\left\{\frac{1}{17}, -\frac{1}{17}, -1\right\}$ | $\{-0.471, -0.067, 0.255\}$ | $\times$ |
| [120, 38] | $\left\{\frac{1}{18}, -\frac{1}{18}, -1\right\}$ | $\{-0.489, -0.074, 0.280\}$ | $\times$ |
| [120, 40] | $\left\{\frac{1}{19}, -\frac{1}{19}, -1\right\}$ | $\{-0.503, -0.078, 0.300\}$ | $\times$ |
| [120, 42] | $\left\{\frac{2}{5}, -\frac{1}{20}, -\frac{1}{8}\right\}$ [93], $\left\{-\frac{2}{5}, -\frac{1}{20}, \frac{1}{8}\right\}$ [93] | $\{-0.515, -0.082, 0.317\}$ | $\times$ |
| | $\left\{\frac{1}{20}, -\frac{1}{20}, -1\right\}$ | | $\times$ |
| [120, 44] | $\left\{\frac{1}{21}, -\frac{1}{21}, -1\right\}$ | $\{-0.525, -0.084, 0.331\}$ | $\times$ |
| [120, 46] | $\left\{\frac{1}{22}, -\frac{1}{22}, -1\right\}$ | $\{-0.534, -0.087, 0.343\}$ | $\left\{\frac{1}{22}, 0.385, 0.046\right\}$, $\left\{-\frac{1}{22}, 0.353, -0.046\right\}$ |

| [$N, N_1$] | Feasible | 3-inner product set | |
|---|---|---|---|
| | | Numerical | |
| | | constant weight | non-constant weight |
| [120, 48] | $\left\{\frac{1}{23}, -\frac{1}{23}, -1\right\}$ | $\{-0.542, -0.088, 0.354\}$ | $\left\{\frac{1}{23}, 0.399, 0.044\right\}$ $\left\{-\frac{1}{23}, 0.365, -0.044\right\},$ $\{-1, 0.308 + 0.007i, -1.001 + 0.007i\},$ |
| [120, 50] | $\left\{\frac{1}{24}, -\frac{1}{24}, -1\right\}$ | $\{-0.549, -0.090, 0.363\}$ | $\{-1, 0.317 - 0.015i, -1.000 - 0.016i\}$ |
| [120, 52] | $\left\{\frac{1}{5}, -\frac{1}{5}, -1\right\}$ [93] | $\{-0.556, -0.091, 0.372\}$ | $\{-1, 1.075, -1.049 - 1.033i\}$ |
| | $\left\{\frac{1}{25}, -\frac{1}{25}, -1\right\}$ | | $\{-1, 1.075, -1.049 - 1.033i\}$ |
| [120, 54] | $\left\{\frac{1}{26}, -\frac{1}{26}, -1\right\}$ | $\{-0.562, -0.092, 0.379\}$ | $\{-1, 0.087 + 0.313i, -1.184\}$ |
| [120, 56] | $\left\{\frac{7}{15}, -\frac{1}{15}, -\frac{3}{5}\right\}$ [93], $\left\{\frac{1}{15}, -\frac{1}{6}, \frac{3}{10}\right\}$ [93], $\left\{-\frac{1}{5}, \frac{1}{10}, -\frac{1}{10}\right\}$ [93] | $\{-0.567, -0.093, 0.386\}$ | × |
| | $\left\{\frac{1}{9}, -\frac{1}{9}, -1\right\}$ [94], $\left\{-\frac{1}{3}, \frac{1}{3}, -1\right\}$ [94] | | × |
| | $\left\{\frac{1}{27}, -\frac{1}{27}, -1\right\}$ | | × |
| [120, 58] | $\left\{\frac{1}{28}, -\frac{1}{28}, -1\right\}$ | $\{-0.571, -0.094, 0.392\}$ | × |
| [120, 60] | $\left\{\frac{1}{29}, -\frac{1}{29}, -1\right\}$ | $\times (R = 1)$ | $\{0.034, 0.452, 0.035\},$ $\{-0.034, 0.418, -0.035\}$ |

We see that none of numerical 3-inner product sets coincide with the feasible ones. We note that in case $N_1 = 42$ (constant weight) our assumption that $w(x) = 1$, for all $x \in X$, implies the square of second radius $R = 1$ which is impossible.